\documentclass[english,a4paper,10pt]{article}
\usepackage{theorem}
\usepackage{amssymb,amsmath}
\usepackage{latexsym}
\usepackage{babel}
\usepackage{color}
\usepackage[latin1]{inputenc}
\usepackage{color,graphicx,hyperref,array,fleqn}

\newtheorem{theorem}{Theorem}[section]
\newtheorem{proposition}[theorem]{Proposition}

\newtheorem{lemma}[theorem]{Lemma}
\theorembodyfont{\normalfont\mdseries\upshape}

\newtheorem{example}[theorem]{Example}
\newtheorem{definition}[theorem]{Definition}

\title{\bf Revisiting the Hahn--Banach Theorem and Nonlinear Infinite Programming}
\author{P. Montiel L\'opez$^{\hbox{(a)}}$ and M. Ruiz Gal\'an$^{\hbox{(b)}}$}

\begin{document}

\maketitle

\centerline{$^{\hbox{(a)}}$University of Granada,
Centro de Estudios Superiores La Inmaculada,} 

\centerline{Department of Sciences,
c/ Joaquina Eguaras, 114, 18013 Granada (Spain),}

\centerline{e-mail: pablomontiel@eulainmaculada.com}

\medskip

\centerline{$^{\hbox{(b)}}$University of Granada,
E.T.S. Ingenier\'{\i}a de Edificaci\'on,}

\centerline{ Department of Applied Mathematics,
c/ Severo Ochoa s/n, 18071 Granada (Spain),}

\centerline{e-mail:
mruizg@ugr.es}

\begin{abstract}
\noindent The aim of this paper is to state a sharp version of the K\"onig supremum theorem, an equivalent reformulation of the Hahn--Banach theorem. We apply it to derive statements of the Lagrange multipliers, Karush--Kuhn--Tucker and Fritz John types, for nonlinear infinite programs. We also show that a weak concept of convexity coming from minimax theory, infsup-convexity, is the adequate one for this kind of results.
 \end{abstract}



\vspace{0.5cm}

\section{Introduction}\label{sect1}

Without a doubt, the Hahn--Banach theorem is not only a central result in functional and convex analysis, but also provides endless applications 
in many other fields, even outside of mathematics. One of its  powerful equivalent reformulations is the so-called \textit{Mazur--Orlicz} theorem (see \cite[Th\'eor\`{e}me 2.41]{maz-orl}, \cite[Theorem, p. 365]{pta} and \cite[Satz, p. 482]{ko02}, \cite[Theorem 28]{sim2}, and its generalizations \cite[Satz, p. 482 and Zusatz, p.483]{ko02}, \cite[Theorem 1.1]{ko0}, \cite[Theorem 2.9]{sim1}, \cite[Theorem 2]{liu}, \cite[Theorem 12]{grz-prz-urb}, \cite[Theorem 3.1]{din-mo}), \cite[Theorem 3.5]{sun} and \cite[Theorem 3.5 and Theorem 6.1]{sim-2}), which allows one to find a linear functional dominated by another sublinear functional, and states in addition a control of the infimum of both functionals on a given convex set. Such a control is not trivial, and generates numerous applications in functional analysis, minimax theory, variational analysis, monotone multifunctions theory or optimization. For instance, one can check \cite{sim1,liu,ga-ru,din-mo}. Along these same lines one can consider the \textit{K\"onig supremum theorem}  
(\cite[Erweiterter Maximumssatz p. 501]{ko01}), which is an extension result in a space of bounded functions, although it is again equivalent to the Hahn--Banach theorem (the details can be found in Proposition \ref{pr:equiv}).

The main result in this paper, Theorem \ref{th:main}, establishes as a consequence of the Mazur--Orlicz theorem, a generalization of the K\"onig supremum theorem in terms of a not very restrictive kind of convexity (\textit{infsup-convexity}, see Definition \ref{de:infsup} below). This notion of convexity arises in minimax theory: see 
\cite[Definition 2.11]{ste}, \cite[p.
653]{ka-ko} and \cite[Definition 2.1]{rui1}. In this setting, infsup-convexity is the adequate type of convexity to state some general
characterizations of the minimax inequality (\cite[Corollary
3.12]{rui2}). Finally, from such a generalization we deduce several theorems for nonlinear infinite programs --Lagrange multipliers, Karush--Kuhn--Tucker, Fritz John--, extending those in the finite case in \cite{rui3,lo-ru}.

\vspace{0.5cm}

\section{Preliminaries}\label{sect2}

Let us begin by evoking the Mazur--Orlicz theorem. Recall that a real-valued functional on a real linear space is \textit{sublinear} if it is subadditive and positively homogeneous.

\bigskip

\begin{theorem}[Mazur--Orlicz]\label{th:mazurorlicz}
Suppose that $E$ is a real vector space, $C$ is a nonempty and convex subset of $E$, and that $S: E \longrightarrow \mathbb{R}$ is a sublinear
functional. Then, there exists a linear functional $L: E
\longrightarrow \mathbb{R}$ such that
\[
x \in E \ \Rightarrow \ L(x) \le S(x)
\]
and
\[
\inf_{x \in C} L(x) = \inf_{x \in C} S(x).
\]
\end{theorem}

\bigskip

Infinite values are allowed in the equality.

Let us also mention the K\"onig supremum theorem, which we will generalize in Section \ref{sect3}. To this end, if $\Lambda$ is a nonempty set, $\Delta_\Lambda$ stands for the subset of the topological dual $\ell^\infty (\Lambda)^*$ of the real Banach space $\ell^\infty (\Lambda)$ (usual sup-norm) of all the bounded real-valued functions defined on $\Lambda$ 
\[
\Delta_\Lambda := \{ \Phi \in \ell^\infty (\Lambda)^* : \ \Phi \le \sup_{\Lambda} \},
\]
that is, $\Phi \in \ell^\infty (\Lambda)^*$ belonging to $\Delta_\Lambda$ means that
\[
\varphi \in \ell^\infty (\Lambda) \ \Rightarrow \ \Phi (\varphi) \le \sup_{\lambda \in \Lambda} \varphi (\lambda).
\]
It is very easy to prove that $\Phi \in \Delta_\Lambda$ if, and only if, $\Phi$ is positive, i.e., $\varphi \ge 0 $ on $\Lambda$ implies $\Phi (\varphi) \ge 0$, and $\Phi (\mathbf{1})=1$, where $\mathbf{1} \in \ell^\infty (\Lambda)$ is the constant function 1. In particular, when $\Lambda$ is a nonempty finite set, $\ell^\infty(\Lambda)$ is of the form $\mathbb{R}^N$ for some $N \in \mathbb{N}$, and $\Delta_\Lambda$ is the probability simplex
\[
\Delta_N:= \{ (t_1,\dots,t_N) \in \mathbb{R}^N : \ t_1,\dots,t_N \ge 0  \hbox{ and } \sum_{j=1}^N t_j =1 \}.
\] 
As we will check in the proof of Proposition \ref{pr:equiv}, the elements in $\Delta_\Lambda$ act as extension functionals. 

Specifically, the K\"onig supremum theorem reads as follows:

\bigskip

\begin{theorem}[K\"onig]\label{th:konig}
Let $E$ be a real linear space, $\Lambda$ be a nonempty set, $L: E \longrightarrow \mathbb{R}$ be a linear functional and assume that, for each $\lambda \in \Lambda$, $S_\lambda : E \longrightarrow \mathbb{R}$ is a sublinear functional in such a way that
\[
x \in E \ \Rightarrow \ (S_\lambda(x))_{\lambda \in \Lambda} \in \ell^\infty (\Lambda)
\]  
and
\[
x \in E \ \Rightarrow \ L(x) \le \sup_{\lambda \in \Lambda} S_\lambda (x).
\]
Then, there exists $\Phi \in \Delta_\Lambda$ satisfying
\[
x \in E \ \Rightarrow \ L(x) \le \Phi ((S_\lambda(x))_{\lambda \in \Lambda}).
\]
\end{theorem}

\bigskip

Even in the finite case, this result implies a wide variety of applications: see \cite{ko01,ko0,lo-ru}.

Before stating a generelization of the preceding result, which will turn out to be sharp in terms of the functions under consideration, let us observe that not only is it a consequence of the Hahn--Banach theorem (see the proof of K\"onig in \cite[Erweiterter Maximumssatz p. 501]{ko01} from a variant of the Mazur--Orlicz theorem), but also an equivalent reformulation. Now we exactly prove that the validity of Theorem \ref{th:konig} implies (and therefore is equivalent to) that of the norm preserving extension version of the Hahn--Banach theorem.

\bigskip

\begin{proposition}\label{pr:equiv}
The K\"onig supremum theorem implies the Hahn--Banach theorem.
\end{proposition}

\noindent \textsc{Proof.} Suppose that $E$ is a real normed space, $F$ is a vector subspace of $E$, and that $y_0^*: F \longrightarrow \mathbb{R}$ is a continuous and linear functional. We are going to find a continuous and linear functional $x_0^*:  E \longrightarrow \mathbb{R}$ with
\[
{x_0^*}_{|F}=y_0^*
\]
and
\[
\| x_0^*\| = \|y_0^*\|
\]
(usual dual norms). 
Assume that $\|y_0^*\|=1$; then it suffices to consider the linear space $F$, the linear functional $y_0^*$, the index set
\[
\Lambda:=\{x^*\in E^*: \ \|x^*\| \le 1 \},
\]
and, for all $(x^*,y)\in \Lambda \times F$,
\[
S_{x^*}(y):=x^*(y),
\]
which clearly satisfy the assumptions in K\"onig's supremum theorem. Then, there exists $\Phi \in \Delta_\Lambda$ such that
\begin{equation}\label{eq:konig}
y \in F \ \Rightarrow \ y_0^*(y) \le \Phi((x^*(y))_{x^* \in \Lambda}).
\end{equation}
Now we can construct the required $x_0^*\in E^*$. Let $\rho : E \longrightarrow \ell^\infty (\Lambda)$ be the linear operator assigning to each $x \in E$ the function $\rho (x): \Lambda \longrightarrow \mathbb{R}$ given for all $x^* \in \Lambda $ by
\[
\rho(x)(x^*):=x^*(x).
\]
It is clear that $\rho$ is well defined and moreover is an linear isomorphism from $E$ into $\ell^\infty (\Lambda)$. Then we can define $x_0^*$ at each $x \in E$ as
\[
x_0^*(x):=\Phi (\rho(x)).
\]
This functional is obviously linear and in addition is continuous, since given $x \in E$,
\[
\begin{array}{rl}
x_0^*(x) & = \Phi (\rho(x)) \\
           & \le \displaystyle \sup_{x^* \in \Lambda}  \rho(x)(x^*), \qquad \hbox{(since } \Phi \in \Delta_\Lambda \hbox{)} \\ 
		   & \le \| x \|,
\end{array}
\]
and thus, in particular,
\[
\| x_0^* \| \le 1.
\]
In view of this inequality, it only remains to show that $x_0^*$ extends $y_0^*$ to $E$. But that is true, because for $y \in F$ it follows from \eqref{eq:konig} that
\[
\begin{array}{rl}
y_0^*(y) & \le \Phi ((x^*(y))_{x^* \in \Lambda}) \\
		 & = \Phi (\rho (y)) \\
         & = x_0^*(y).
\end{array}
\]
\hfill${\Box}$

\bigskip

In the next section we will establish a sharp version of the K\"onig supremum theorem in terms of the following weak notion of convexity,  useful in minimax theory as mentioned in the Introduction: 

\bigskip

\begin{definition}\label{de:infsup}
Let $X$ and $\Lambda$ be nonempty sets and for each $\lambda \in \Lambda$ let  $f_\lambda : X \longrightarrow \mathbb{R}$ be a function. The family $(f_\lambda)_{\lambda \in \Lambda}$ is said to be
\textit{infsup-convex} on $X$ provided that
\[
\left. 
\begin{array}{c}
m \ge 1, \ \mathbf{t} \in \Delta_m \\ 
x_1,\dots,x_m \in X
\end{array}
\right\} \ \Rightarrow \
\inf_{x \in X} \sup_{\lambda \in \Lambda} f_\lambda(x) \le \sup_{\lambda \in \Lambda} 
\sum_{j=1}^m t_j f_\lambda(x_j).
\]

\end{definition}

\bigskip

Infsup-convexity not only properly extends the notion of convexity of a family functions, but also that of convexlikeness for a family
of functions, due to K. Fan (\cite[p. 42]{fan}).

\vspace{0.5cm}

\section{A Generalized Version of K\"onig's Supremum Theorem}\label{sect3}

Now we focus on stating a general K\"onig supremum theorem, by replacing the vector space with a set, and the linear functional and the family of sublinear functionals with a suitable infsup-convex family of functions, so no linear structure is required. Furthermore, we prove that such a result is sharp.

\bigskip

\begin{theorem}\label{th:main}
Let $X$ and $\Lambda$ be nonempty sets, $f: X \longrightarrow \mathbb{R}$ be a function and $(f_\lambda)_{\lambda \in \Lambda}$ be a family of real valued functions defined on $X$ such that the family $(f_\lambda-f)_{\lambda \in \Lambda}$ is infsup-convex on $X$. Assume in addition that
\[
x \in X \ \Rightarrow \ (f_\lambda (x))_{\lambda \in \Lambda} \in \ell^\infty (\Lambda)
\]
and
\[
x \in X \ \Rightarrow \ f(x) \le \sup_{\lambda \in \Lambda} f_\lambda(x).
\]
Then there exists $\Phi \in \Delta_\Lambda$ such that
\[
x \in X \ \Rightarrow \ f(x) \le \Phi ((f_\lambda)_{\lambda \in \Lambda}).
\]
\end{theorem}

\noindent \textsc{Proof.} Apply the Mazur--Orlicz theorem, Theorem \ref{th:mazurorlicz}, to the real vector space $\ell^\infty (\Lambda)$, its nonempty convex subset
\[
C:=\mathrm{conv}\{(f_\lambda(x)-f(x))_{\lambda \in \Lambda}: \ x \in X \},
\]
and the sublinear functional $S: \ell^\infty(\Lambda) \longrightarrow \mathbb{R}$ given for each $\varphi \in \ell^\infty(\Lambda)$ by
\[
S(\varphi):= \sup_{\lambda \in \Lambda} \varphi(\lambda).
\]
Then, there exists $\Phi \in \Delta_\Lambda$ such that
\[
\inf_{\varphi\in C} \Phi (\varphi) = \inf_{\varphi\in C} S(\varphi).
\]
To conclude, let us observe, on the one hand, that
\[
\inf_{\varphi \in C} \Phi(\varphi) = \inf_{x\in X} \Phi((f_\lambda(x)-f(x))_{\lambda \in \Lambda}),
\]
and on the other hand, that the infsup-convexity of the family $(f_\lambda-f)_{\lambda \in \Lambda}$ on $X$ and the assumption $f \le S((f_\lambda)_{\lambda \in \Lambda})$ yield
\[
\begin{array}{rl}
\displaystyle \inf_{\varphi\in C} S(\varphi) & =  \displaystyle\inf_{\substack{m\ge 1 , \mathbf{t} \in \Delta_m \\  x_1,\dots,x_m \in X}} \sup_{\lambda \in \Lambda}  \sum_{j=1}^m t_j(f_\lambda(x_j)-f(x_j)) \\                                                   & \ge \displaystyle \inf_{x \in X} \sup_{\lambda \in \Lambda} (f_\lambda(x)-f(x)) \\
                                                   & \ge 0.
\end{array}
\]
Therefore, for some $\Phi \in \Delta_\Lambda$,
\[
0 \le  \inf_{x\in X} \Phi((f_\lambda(x)-f(x))_{\lambda \in \Lambda}),
\]
and taking into account that $\Phi$ is linear and $\Phi (\mathbf{1})=1$, because $\Phi \in \Delta_\Lambda$, then
\[
x \in X \ \Rightarrow \ f(x) \le \Phi ((f_\lambda)_{\lambda \in \Lambda}),
\]
as announced. \hfill${\Box}$

\bigskip

The equivalence of Theorem \ref{th:main} and the Hahn--Banach theorem follows from that of the Hahn--Banach theorem and the Mazur--Orlicz theorem and from Proposition \ref{pr:equiv}.

Since infsup-convexity is a concept invariant by adding a constant, then Theorem \ref{th:main} can be equivalently reformulated as follows: assume that $\alpha \in \mathbb{R}$, $X$ and $\Lambda$ are nonempty sets, $f: X \longrightarrow \mathbb{R}$ is a function, and that $(f_\lambda)_{\lambda \in \Lambda}$ is a family of real valued functions defined on $X$ such that the family $(f_\lambda-f)_{\lambda \in \Lambda}$ is infsup-convex on $X$ and for all $x \in X$, $(f_\lambda (x))_{\lambda \in \Lambda} \in \ell^\infty (\Lambda)$. If in addition
\[
x \in X \ \Rightarrow \ f(x) + \alpha  \le \sup_{\lambda \in \Lambda} f_\lambda(x),
\]
then, there exists $\Phi \in \Delta_\Lambda$ such that
\[
x \in X \ \Rightarrow \ f(x)+\alpha \le \Phi ((f_\lambda)_{\lambda \in \Lambda}).
\]
Surprisingly, the converse is also true, as we state in the following result, which is the above-mentioned sharpness of Theorem \ref{th:main}: 

\bigskip

\begin{theorem}\label{th:characterization}
Suppose that $X$ and $\Lambda$ are nonempty sets, $f: X \longrightarrow \mathbb{R}$ is a function, and that $(f_\lambda)_{\lambda \in \Lambda}$ is a family of real valued functions defined on $X$ such that, for each $x \in X$, $(f_\lambda(x))_{\lambda \in \Lambda} \in \ell^\infty (\Lambda)$. Then, the family $(f_\lambda-f)_{\lambda \in \Lambda}$ is infsup-convex on $X$ if, and only if, for all $\alpha \in \mathbb{R}$ satisfying
\[
x \in X \ \Rightarrow \ f(x) + \alpha \le \sup_{\lambda \in \Lambda} f_\lambda (x),
\] 
there exists $\Phi \in \Delta_\Lambda$ such that
\[
x \in X \ \Rightarrow \ f(x)+\alpha \le \Phi ((f_\lambda(x))_{\lambda \in \Lambda}).
\]
\end{theorem}

\noindent \textsc{Proof.} According to the preceding argument, we focus on proving the sufficiency. Hence, let $m\ge
1$, $\mathbf{t} \in \Delta_m$ and $x_1,\dots,x_m \in X$. Let $\alpha:= \inf_{x \in X} \sup_{\lambda \in \Lambda} (f_\lambda(x)-f(x))$, which can be assumed finite without any loss of generality. Since for all $x \in X$,
\[
\begin{array}{rl}
f(x)+\alpha & = \displaystyle f(x)+\inf_{x \in X} \sup_{\lambda \in \Lambda} (f_\lambda(x)-f(x)) \\
            & \le \displaystyle \sup_{\lambda \in \Lambda} f_\lambda (x),
\end{array}
\]
in view of our assumption, we arrive at
\[
x \in X \ \Rightarrow \ f(x)+\alpha \le \Phi ((f_\lambda(x))_{\lambda \in \Lambda})
\]
for some $\Phi \in \Delta_\Lambda$, and therefore, 
\[
\begin{array}{rl}
  \alpha  & \displaystyle \le \inf_{x \in X} \Phi ((f_\lambda(x)-f(x))_{\lambda \in \Lambda}) \qquad \hbox{(since } \Phi (\mathbf{1})=1  \hbox{)} \\
   & \displaystyle \le \min_{j=1,\dots,m} \Phi ((f_\lambda(x_j)-f(x_j))_{\lambda \in \Lambda})  \\
   & \displaystyle \le \sum_{j=1}^m t_j \Phi((f_\lambda(x_j)-f(x_j))_{\lambda \in \Lambda})  \\
   & \displaystyle = \Phi \left(\sum_{j=1}^m t_j (f_\lambda(x_j)-f(x_j))_{\lambda \in \Lambda} \right)  \\   
   & \displaystyle \le \sup_{\lambda \in \Lambda}\sum_{j=1}^m t_j (f_\lambda(x_j)-f(x_j)) \qquad \hbox{(because } \Phi \le \displaystyle \sup_\Lambda \hbox{).}
\end{array}
\]
The arbitrariness of $m\ge
1$, $\mathbf{t} \in \Delta_m$ and $x_1,\dots,x_m \in X$ yields the announced infsup-convexity.
\hfill${\Box}$

\bigskip

Let us point out that Theorem \ref{th:main} and Theorem \ref{th:characterization} were proven for $\Lambda$ finite in \cite[Theorem 2.3]{lo-ru} and \cite[Theorem 2.4]{lo-ru}, respectively.

\vspace{0.5cm}

\section{Consequences in Infinite Programming}\label{sect4}

Assume that $X$ and $\Lambda $ are nonempty sets, $f: X \longrightarrow \mathbb{R}$,  
$(f_\lambda)_{\lambda \in \Lambda}$ is a family of real valued functions on $X$ satisfying
\[
x \in X \ \Rightarrow \ (f_\lambda (x))_{\lambda \in \Lambda} \in \ell^\infty (\Lambda),
\]  
and that the set 
\[
X_0:=\left\{ x \in X : \ \sup_{\lambda \in \Lambda} f_\lambda(x) \le 0 \right\}
\]
is nonempty. Let us consider the nonlinear infinite program
\begin{equation}\label{eq:nlp}
\inf_{x \in X_0} f(x).
\end{equation}
If we denote by $\ell^\infty (\Lambda)^*_+$ the cone of the positive functionals in $\ell^\infty (\Lambda)^*$, then the associated \textit{Lagrangian} $\mathbf{L} : X \times \ell^\infty (\Lambda)_+^* \longrightarrow \mathbb{R}$ is defined at each $(x,\Phi) \in X \times \ell^\infty (\Lambda)^*_+$ as
\[
\mathbf{L}(x,\Phi):=f(x)+\Phi\left( (f_\lambda(x))_{\lambda \in \Lambda} \right).
\] 
In addition, $(x^0,\Phi_0) \in X \times \ell^\infty (\Lambda)^*_+$ is said to be a \textit{saddle point} of $\mathbf{L}$ provided that
\[
(x,\Phi) \in X \times \ell^\infty (\Lambda)^*_+ \ \Rightarrow \ \mathbf{L}(x^0,\Phi) \le \mathbf{L}(x^0,\Phi_0) \le \mathbf{L}(x,\Phi_0).
\]
In such a case, $\Phi_0$ is a \textit{Lagrange multiplier} for $\mathbf{L}$.
It is an elementary fact that $x^0 \in X$ is an optimal solution of the nonlinear problem \eqref{eq:nlp} provided there exists $\Phi_0 \in \ell^\infty (\Lambda)^*_+$ such that $(x^0, \Phi_0)$ is a saddle point for the Lagrangian. Now we go the other way, by proving that the infsup-convexity of a certain family of functions is exactly the assumption required for deriving --under a natural Slater condition-- the existence of a Lagrange multiplier $\Phi_0 \in \ell^\infty (\Lambda)^*_+$ from that of an optimal solution $x^0 \in X$. It is a Lagrange multiplier type result: see the classical theorem \cite{uza,roc} and its extensions in convexlike and quasiconvex contexts \cite{arr-ent,hay-kom,ill-kas,suz-kur,fan-luo-wan}.

Before it, an easy technical result:

\bigskip

\begin{lemma}\label{le:easy}
Let $X$ and $\Lambda$ be nonempty sets, $(f_\lambda)_{\lambda \in 
\Lambda}$ be a family of real valued functions defined on $X$ such that
\[
x \in X \ \Rightarrow \ (f_\lambda(x))_{\lambda \in \Lambda} \in \ell^\infty (\Lambda),
\] 
and suppose that the set $X_0:=\left\{ x \in X : \ \sup_{\lambda \in \Lambda} f_\lambda(x) \le 0 \right\}$ is nonempty. If $x^0 \in X$ is a solution of the nonlinear program \eqref{eq:nlp}, then
\[
\inf_{x \in X} \max \left\{\sup_{\lambda \in \Lambda} f_\lambda (x), f(x)-f(x^0) \right\}=0.
\]
\end{lemma}

\noindent \textsc{Proof.} We are assuming that 
\[
f(x^0)=\inf_{x \in {X_0}} f(x),
\]
hence
\[
0 \le \inf_{x \in {X_0}} \max \left\{\sup_{\lambda \in \Lambda} f_\lambda (x), f(x)-f(x^0) \right\}.
\]
But for all $x \in X \backslash X_0$ there exists $\lambda \in \Lambda$ with $0 < f_\lambda (x)$, so
\[
0 \le \inf_{x \in X \backslash {X_0}} \max \left\{\sup_{\lambda \in \Lambda} f_\lambda (x), f(x)-f(x^0) \right\}.
\]
According to these two previous inequalities, we arrive at
\[
0 \le \inf_{x \in X} \max \left\{\sup_{\lambda \in \Lambda} f_\lambda (x), f(x)-f(x^0) \right\}
\]
and, since
\[
0 = \max \left\{\sup_{\lambda \in \Lambda} f_\lambda (x^0), f(x^0)-f(x^0) \right\},
\]
we have concluded the proof.
\hfill${\Box}$

\bigskip
  
Now we are in a position to state the aforementioned relationship between optimal solutions and saddle points (or Lagrange multipliers). This is our main statement on nonlinear infinite programming.

\bigskip

\begin{theorem}\label{th:lagrange}
Suppose that $X$ and $\Lambda$ are nonempty sets, $f : X \longrightarrow \mathbb{R}$ and that $(f_\lambda)_{\lambda \in \Lambda}$ is a family of real valued functions defined on $X$ such that
\[
x \in X \ \Rightarrow \ (f_\lambda (x))_{\lambda \in \Lambda} \in \ell^\infty (\Lambda),
\]
and the feasible set $X_0:=\left\{ x \in X : \ \sup_{\lambda \in \Lambda} f_\lambda(x) \le 0 \right\}$ is nonempty. Let us also assume that $x^0\in X$ is an optimal solution of the nonlinear problem \eqref{eq:nlp} and that the following Slater condition is fulfilled:
\[
\hbox{there exists } x^1 \in X \hbox{ such that } \sup_{\lambda \in \Lambda} f_\lambda (x^1) <0.
\]
Then there exists $\Phi_0 \in \ell^\infty (\Lambda)^*_+$ such that $(x^0,\Phi_0)$ is a saddle point for $\mathbf{L}$ if, and only if, the family $(f_\lambda)_{\lambda \in \Lambda} \cup (f-f(x^0))$ is infsup-convex on $X$.
\end{theorem}

\noindent \textsc{Proof.} We first assume that for some $\Phi_0 \in \ell^\infty (\Lambda)^*_+$, $(x^0,\Phi_0)$ is a saddle point of $\mathbf{L}$. Then $f(x^0)=\mathbf{L}(x^0,\Phi_0)$ and 
\[
x \in X \ \Rightarrow \ \mathbf{L}(x^0,\Phi_0) \le \mathbf{L}(x,\Phi_0),
\]
i.e.,
\begin{equation}\label{eq:ecuacion}
0 \le \inf_{x \in X} (\Phi_0 ((f_\lambda(x)))_{\lambda \in \Lambda}+ f(x)-f(x^0)).
\end{equation}
Let $\Lambda_0:=\Lambda \cup \{ \mu \}$, where $\mu \notin \Lambda$, and define $\Psi \in \ell^\infty (\Lambda_0)^*=\ell^\infty (\Lambda)^* \times \mathbb{R}$ as
\[
\Psi:=\frac{1}{1+\Phi_0(\mathbf{1})}(\Phi_0,1),
\]
which clearly belongs to $\Delta_{\Lambda_0}$ and, thanks to \eqref{eq:ecuacion} satisfies
\[
0 \le \inf_{x \in X} \Psi((f_\lambda (x))_{\lambda \in \Lambda} \cup (f(x)-f(x^0))).
\]
But, since $x^0$ is an optimal solution for \eqref{eq:nlp}, in view of Lemma \ref{le:easy},
\[
0=\inf_{x \in X} \max \left\{\sup_{\lambda \in \Lambda} f_\lambda (x), f(x)-f(x^0) \right\}.
\] 
Therefore, we have that
\[
\inf_{x \in X} \max \left\{\sup_{\lambda \in \Lambda} f_\lambda (x), f(x)-f(x^0) \right\} \le \inf_{x \in X} \Psi((f_\lambda (x))_{\lambda \in \Lambda} \cup (f(x)-f(x^0)))
\]
and this inequality and the fact that $\Psi \in \Delta_{\Lambda_0}$ easily imply, as in the last part of the proof of Theorem 
\ref{th:characterization}, the infsup-convexity of the family $(f_\lambda)_{\lambda \in \Lambda} \cup (f-f(x^0))$ on $X$.

And conversely, let us suppose that $x^0\in X$ is an optimal solution of the nonlinear program under consideration and that the Slater condition is satisfied. The first assumption and Lemma \ref{le:easy} yield
\[
\inf_{x \in X} \max \left\{\sup_{\lambda \in \Lambda} f_\lambda (x), f(x)-f(x^0) \right\}=0.
\]
Then, Theorem \ref{th:main}, when applied with the function $f: X \longrightarrow \mathbb{R}$, assigning to each $x \in X$ the value
\[
f(x):=0,
\]
and the infsup-convexity on $X$ of the family $(f_\lambda)_{\lambda \in \Lambda} \cup (f-f(x^0))$, provides us with a positive and linear functional $\Phi : \ell^\infty (\Lambda) \longrightarrow \mathbb{R}$ and $\rho \ge 0$  with
\[
\Phi (\mathbf{1})+\rho=1
\]
and
\begin{equation}\label{eq:desig}
x \in X \ \Rightarrow \ 0 \le \Phi ((f_\lambda(x))_{\lambda \in \Lambda})+\rho (f(x)-f(x^0)).
\end{equation}
Let us notice that $\rho >0$, because otherwise $\Phi \in \Delta_\Lambda$ and we would arrive at
\[
\begin{array}{rl}
0 & \le \displaystyle \inf_{x \in X} \Phi ((f_\lambda(x))_{\lambda \in \Lambda}) \\
  &  \\ 
& \le \displaystyle \Phi ((f_\lambda(x^1))_{\lambda \in \Lambda}) \\
   & \\
& \displaystyle \le \sup_{\lambda \in \Lambda} f_\lambda (x^1) \\
&  \\
& < 0,
\end{array}
\]
a contradiction. So $\rho >0$ and we take $\Phi_0:=\Phi / \rho \in \ell^\infty (\Lambda)^*_+$. Then, according to \eqref{eq:desig} we have that
\[
x \in X \ \Rightarrow \ f(x^0) \le f(x) + \Phi_0((f_\lambda(x))_{\lambda \in \Lambda}).
\]
Then, given $x \in X$,
\[
\begin{array}{rl}
\mathbf{L}(x^0,\Phi_0) & =  f(x^0)+  \Phi_0((f_\lambda(x^0))_{\lambda \in \Lambda}) \\
    &  \\
    & \le f(x^0) \\
    &  \\
    & \le  f(x) + \Phi_0((f_\lambda(x))_{\lambda \in \Lambda}) \\
    &  \\
    & = \mathbf{L}(x, \Phi_0).
\end{array} 
\] 
Finally, given $\Upsilon \in \ell^\infty (\Lambda)^*_+$ it follows that
\[
\Upsilon ((f_\lambda(x^0)_{\lambda \in \Lambda})) \le 0,
\]
so
\[
\mathbf{L}(x^0,\Upsilon) \le \mathbf{L}(x^0,\Phi_0)
\]
and we have completed the proof. 
\hfill${\Box}$

\bigskip

We would hope that the Slater condition in Theorem \ref{th:lagrange} could be replaced with this weaker one:
\begin{equation}\label{eq:contral}
\hbox{there exists } x^1 \in X \hbox{ such that } \lambda \in \Lambda \ \Rightarrow \ f_\lambda (x^1)<0,
\end{equation}
but it is false in general:

\bigskip

\begin{example}
Let $X:=\mathbb{R}$, $\Lambda:=\mathbb{N}$, $f: X \longrightarrow \mathbb{R}$ be the function assigning to each $x \in X$
\[
f(x):=x, 
\]
and for all $n \in \mathbb{N}$, let $f_n: \mathbb{R} \longrightarrow \mathbb{R}$ be the function defined at each $x \in \mathbb{R}$ by
\[
f_n(x)=-\frac{x^3}{n}.
\]
Then, the feasible set is 
\[
X_0=\mathbb{R}_+.
\]
The Slater condition fails, for if there exists $x^1 \in \mathbb{R}$ with
\[
\sup_{n \in \mathbb{N}} f_n(x^1)<0,
\]
then, in particular, $x^1 \in X_0=\mathbb{R}_+$, but for such an $x^1$,
\[
\sup_{n \in \mathbb{N}}f_n(x^1)=0,
\]
which contradicts the Slater condition. However, clearly any $x^1 \in X_0 \backslash \{0\}$ fulfils \eqref{eq:contral}. To conclude, let us show that, despite the fact that the nonlinear program 
\[
\inf_{x \in {X_0}} f(x)
\]
admits the optimal solution $x^0=0$, there exists no a corresponding Lagrange multiplier $\Phi_0 \in \ell^\infty (\mathbb{N})^*_+$ for the associated Lagrangian $\mathbf{L}$. To this end, in order to check the hypotheses in Theorem \ref{th:lagrange} (obviously except the Slater condition), we prove the unique non trivial fact that the family $(f_n)_{n \in \mathbb{N}} \cup f$ is infsup-convex on $\mathbb{R}$ ($f(x^0)=0$), that is,
\[
\inf_{x \in \mathbb{R}} \left( \sup_{n \in \mathbb{N}} f_n(x) \vee f(x) \right) \le \sup_{n \in \mathbb{N}}\sum_{j=1}^m t_j f_n(x_j) \vee \sum_{j=1}^m t_j f(x_j),
\] 
whenever $m \ge 1$, $\mathbf{t} \in \Delta_m$ and $x_1,\dots,x_m \in \mathbb{R}$.
By Lemma \ref{le:easy} we know for the left-hand side that
\[
\inf_{x \in \mathbb{R}} \left( \sup_{n \in \mathbb{N}} f_n(x) \vee f(x) \right) =0,
\]
so we have to show that
\[
0 \le \sup_{n \in \mathbb{N}}\sum_{j=1}^m t_j f_n(x_j) \vee \sum_{j=1}^m t_j f(x_j),
\]
i.e.,
\[
0 \le \sup_{n \in \mathbb{N}}\left( \sum_{j=1}^m t_j x_j^3\right) \left( -\frac{1}{n}\right) \vee \left( \sum_{j=1}^m t_j x_j \right).
\]
But this inequality is clearly satisfied, because if 
\[
0 \le \sum_{j=1}^m t_jx_j^3 
\]
then
\[
\sup_{n \in \mathbb{N}}\left( \sum_{j=1}^m t_j x_j^3\right) \left( -\frac{1}{n}\right)=0
\]
and the inequality holds, while if
\[
\sum_{j=1}^m t_jx_j^3  < 0
\]
then
\[
0 \le \sup_{n \in \mathbb{N}}\left( \sum_{j=1}^m t_j x_j^3\right) \left( -\frac{1}{n}\right),
\]
which implies the validity of the inequality and thus the above mentioned infsup-convexity.

We finish by arguing by contradiction, so let us assume that the nonlinear problem under consideration admits a Lagrange multiplier $\Phi_0 \in {\ell^\infty}(\mathbb{N})^*_+$. Since, in particular, $\Phi_0 \in {\ell^\infty}(\mathbb{N})^*$, making use of the Dixmier decomposition of $\ell^\infty (\mathbb{N})^*$, there exist $y_0 \in \ell^1(\mathbb{N})$ and $\varphi_0 \in c_0(\mathbb{N})^\perp$ such that 
\[
\Phi_0=y_0+\varphi_0,
\]
where $c_0(\mathbb{N})$ is the closed linear subspace of $\ell^\infty (\mathbb{N})$ of those null sequences and $c_0(\mathbb{N})^\perp$ is its annihilator. Observe that for all $n \in \mathbb{N}$, $y_0(n) \ge 0$, because $\Phi_0 \in {\ell^\infty}(\mathbb{N})^*_+$ and $e_n=(0,\dots,0,\overbrace{1}^{n)},0, \dots, 0,\dots) \in \ell^\infty(\mathbb{N})_+$, so
\[
\begin{array}{rl}
0 & \le \Phi_0(e_n)   \\
& = y_0(e_n)+\varphi_0(e_n)  \\
& = y_0(n) \qquad (\varphi_0 \in c_0(\mathbb{N})^\perp, \ e_n \in c_0(\mathbb{N})).
\end{array}
\]
Then, taking into account that we are assuming that $(0,\Phi_0)$ is a saddle point for the Lagrangian, in particular there holds for all $x \in \mathbb{R}$
\[
f(x^0)+\Phi_0 ((f_n(x^0))_{n\in \mathbb{N}}) \le f(x) +\Phi_0 ((f_n(x))_{n\in \mathbb{N}}).
\]
But, for each $x \in \mathbb{R}$ $(f_n(x))_{n\in \mathbb{N}} \in c_0(\mathbb{N})$, so this inequality is nothing more than
\begin{equation}\label{eq:contra}
x \in \mathbb{R} \ \Rightarrow \ 0 \le x - x^3 \sum_{n=1}^\infty \frac{y_0(n)}{n}
\end{equation}
which is absurd, because if $\sum_{n=1}^\infty \frac{y_0(n)}{n}=0$, then it suffices to take $x<0$ in \eqref{eq:contra} to arrive at a contradiction, while if $\sum_{n=1}^\infty \frac{y_0(n)}{n} > 0$, a large enough $x>0$ yields
\[
x - x^3 \sum_{n=1}^\infty \frac{y_0(n)}{n}<0,
\]
once again against the inequality  \eqref{eq:contra}.
\hfill${\Box}$

\end{example}

\bigskip

We finish by deriving from Theorem \ref{th:lagrange} and Theorem \ref{th:main} some Karush--Kuhn--Tucker and Fritz John results (see \cite{kar,joh,kuh-tuc,bre-tre,ito-kun,flo}), respectively, for which the sharp concept of convexity turns out to be once again infsup-convexity. The first of them 
is just an equivalent reformulation of Theorem \ref{th:lagrange}, according to the easy-to-prove fact that the Karush--Kuhn--Tucker conditions below are equivalent to the existence of a Lagrange multiplier:

\bigskip

\begin{theorem}
Assume that $X$ and $\Lambda$ are nonempty sets, $x^0 \in X$, $f: X \longrightarrow \mathbb{R}$ and that $(f_\lambda)_{\lambda \in \Lambda}$ is a family of real valued functions defined on $X$ in such a way that the feasible set $X_0:=\{x \in X: \ \sup_{\lambda \in \Lambda} f_\lambda(x) \le 0 \}$ is nonempty and that
\[
x \in X \ \Rightarrow \ (f_\lambda(x))_{\lambda \in \Lambda} \in \ell^\infty (\Lambda).
\]
If in addition the family $(f_\lambda)_{\lambda \in \Lambda} \cup (f-f(x^0))$ is infsup-convex on $X$ and the Slater condition 
\[
\hbox{there exists } x^1 \in X \hbox{ such that } \sup_{\lambda \in \Lambda} f_\lambda (x^1) <0
\]
is valid, then $x^0$ is an optimal solution for the nonlinear problem \eqref{eq:nlp} if, and only if, $x^0 \in X_0$ and
there exists $\Phi_0 \in \ell^\infty (\Lambda)^*_+$ such that
\[
f+\Phi_0((f_\lambda(\cdot))_{\lambda \in \Lambda}) \hbox{ attains its infimum on } X \hbox{ at } x^0
\]
and
\[
\Phi_0((f_\lambda(x^0))_{\lambda \in \Lambda})=0.
\]
\end{theorem} 

\bigskip

The Fritz John result needs some additional work:

\bigskip

\begin{theorem}
Let $X$ and $\Lambda$ be nonempty sets, $f: X \longrightarrow \mathbb{R}$ and for all $\lambda \in \Lambda$, $f_\lambda : X \longrightarrow \mathbb{R}$ such that the feasible set $X_0:=\{x \in X: \ \sup_{\lambda \in \Lambda} f_\lambda(x) \le 0\}$ is nonempty, $x^0$ is an optimal solution of the infinite program \eqref{eq:nlp}, and
\[
x \in X \ \Rightarrow \ (f_\lambda(x))_{\lambda \in \Lambda} \in \ell^\infty (\Lambda).
\]
Then, there exist $\rho \ge 0$ and $\Phi_0 \in \ell^\infty(\Lambda)^*_+$ with $\rho+\Phi_0(\mathbf{1})=1$ satisfying the following Fritz John conditions
\begin{equation}\label{eq:fj1}
\rho f + \Phi_0((f_\lambda(\cdot))_{\lambda \in \Lambda}) \hbox{ attains its infimum on } X \hbox{ at } x^0
\end{equation}
and
\begin{equation}\label{eq:fj2}
\Phi_0((f_\lambda(x^0))_{\lambda \in \Lambda})=0
\end{equation}
if, and only if, the family $(f_\lambda)_{\lambda \in \Lambda} \cup (f-f(x^0))$ is infsup-convex on $X$.
\end{theorem}

\noindent \textsc{Proof.} Suppose that the family $(f_\lambda)_{\lambda \in \Lambda} \cup (f-f(x^0))$ is infsup-convex on $X$. Apply Lemma \ref{le:easy} to arrive at
\[
\inf_{x \in X} \max \left\{\sup_{\lambda \in \Lambda} f_\lambda (x), f(x)-f(x^0) \right\}=0.
\]
Then, Theorem \ref{th:main} provides us with a $\Phi \in \Delta_{\Lambda_0}$ (same notation as in the proof of Theorem \ref{th:lagrange}) such that
\[
0 \le \Phi ((f_\lambda)_{\lambda \in \Lambda} \cup (f-f(x^0))),
\]
i.e., for some $\rho \ge 0$ and $\Phi_0 \in \ell^\infty (\Lambda)^*_+$ we have that
\[
\Phi_0(\mathbf{1}) \ge 0,
\]
\[
\Phi_0(\mathbf{1})+ \rho=1
\]
and 
\begin{equation}\label{eq:1}
x \in X \ \Rightarrow \ \rho f(x^0) \le \rho f(x)+\Phi_0((f_\lambda(x))_{\lambda \in \Lambda}).
\end{equation}
Condition \eqref{eq:fj2} follows from both the fact that $x^0 \in X_0$ and $\Phi_0 \in \ell^\infty (\Lambda)^*_+$, and inequality \eqref{eq:1} for $x=x^0$. And conditions \eqref{eq:fj2} and \eqref{eq:1} clearly imply \eqref{eq:fj1}.

And conversely, if $\rho \ge 0$ and $\Phi_0 \in \ell^\infty(\Lambda)^*_+$ with $\rho+\Phi_0(\mathbf{1})=1$ fulfil conditions \eqref{eq:fj1} and \eqref{eq:fj2}, then
\[
x \in X \ \Rightarrow \ 0 \le \rho (f(x)-f(x^0))+ \Phi_0 ((f_\lambda (x))_{\lambda \in \Lambda}). 
\]
In particular, $\Psi:=(\Phi_0,\rho) \in \Delta_{\Lambda_0}$ (notation in the proof of Theorem \ref{th:lagrange}) satisfies, according to Lemma \ref{le:easy}, that 
\[
\inf_{x \in X} \max \left\{\sup_{\lambda \in \Lambda} f_\lambda (x), f(x)-f(x^0) \right\} \le \inf_{x \in X} \Psi((f_\lambda (x))_{\lambda \in \Lambda} \cup (f-f(x^0))),
\]
which, as mentioned in the proof of Theorem \ref{th:lagrange}, implies the infsup-convexity of the family $(f_\lambda)_{\lambda \in \Lambda} \cup (f-f(x^0))$ on $X$.
\hfill${\Box}$

\bigskip

\vspace{0.5cm}

\section*{Acknowledgement}

Research partially supported by project MTM2016-80676-P (AEI/FEDER, UE), Junta de Andaluc\'{\i}a Grant FQM359, and Centro de Estudios Superiores La Inmaculada and E.T.S. de Ingenier\'{\i}a de Edificaci\'on (University of Granada, Spain).

\end{document}